\documentclass[11pt]{pjm} 
\usepackage{graphics} 

\title[Convexity for  the figure eight]
{Convexity in the figure eight solution to
the three-body problem} 

\author[T. Fujiwara]{Toshiaki Fujiwara} 
\address{Faculty of General Studies\\
Kitasato University\\
Kitasato 1-15-1, Sagamihara\\
Kanagawa 228-8555, Japan} 
\email{fujiwara@clas.kitasato-u.ac.jp} 
\author[R. Montgomery]{Richard Montgomery} 
\address{Mathematics Dept.\\
UC Santa Cruz\\
Santa Cruz, CA 95064, USA} 
\email{rmont@count.ucsc.edu} 

\thanks{The authors would like to thank  AIM/ARCC for
funding a workshop in celestial  mechanics where the authors
met.  We would also like to 
 express our sincere thanks to
 Alain Chenciner, Hiroshi Fukuda and Hiroshi Ozaki
for discussions.}

\newtheorem{theorem}{Theorem}
\newtheorem{lemma}{Lemma}
\newtheorem{corollary}{Corollary}
\newtheorem{proposition}{Proposition}

\def\R{I\!\! R}

\date{Received date / Revised version date} 

\begin{document} 

\begin{abstract} 

The figure eight is a remarkable solution to the Newtonian
three-body problem in which the three equal masses 
chase each around a planar curve having the qualitative
shape and symmetries of a figure eight.  
Here we prove
that each lobe of this eight is convex.   

\end{abstract} 

\maketitle 

\section{Introduction}\label{intro} 
The figure eight is a recently discovered periodic  solution to the
Newtonian   three-body problem   in which     three 
equal masses
traverse a single closed planar curve   which has the 
form 
of a figure eight (figure \ref{theInterval}). 
See \cite{Moore}, and \cite{CM}.
In particular,
it has one point of self-intersection,
the origin, which divides the eight
into two symmetric parts, its two 
lobes.
In \cite{CM} it was
proved that each lobe is star-shaped.
Here we prove  convexity of the lobes. 

\begin{theorem}\label{theo:convexity}
 Each lobe of the eight solution
is a convex curve.
\end{theorem}

In the final section we describe  how the theorem generalizes
to prove the convexity of eights
for many three body potentials besides Newton's.

A computer proof based on interval arithmetic appears in \cite{KZ}.

\section{Preliminaries.}\label{preliminaries}
We present  a number of  
properties of the eight established in \cite{CM}
and three assertions relating   mechanics and plane
geometry.  The convexity proof relies on these properties and assertions. 

\subsection{Center of Mass.}\label{centerOfMass}
 Write $(q_1 (t), q_2 (t), q_3 (t))$
for the location of the three masses at time 
$t$. 
The $q_i (t)$ are points in the plane.
At each time $t$
we have that 
$q_1 (t) + q_2 (t) + q_3 (t) = 0$.

\subsection{Symmetry.}\label{symmetry}
Write  
$R_y( x,y) = (-x,y)$
for the reflection about the   $y$ axis. 
Then the eight solution enjoys the following symmetries: 
$$
(q_1 (t), q_2 (t), q_3 (t))  =  
(R_y (q_3 (t-T/6)),R_y (q_1(t-T/6)),R_y (q_2(t-T/6)) )$$
$$(q_1 (t), q_2 (t), q_3 (t)) =
(-q_1(-t),-q_3(-t), -q_2(-t)).$$ 
  The right-hand side of these equations 
define transformations
 $s$ and $\sigma$ on the space of all $T$-periodic
loops.  
These transformations generate
an action of the    dihedral  group 
$$D_6 = \{s,\sigma| s^6=1, \sigma^2=1, 
s\sigma=\sigma s^{-1} \},$$ 
the symmetry group  of a regular hexagon, 
which is consequently a symmetry
group of the eight.  

Invariance under   $s^2 \in D_6$
implies that   $(s^2 (q_1, q_2, q_2) )(t) =  (q_1 (t), q_2 (t),
q_3(t))$. 
Setting  $q = q_1$
this last equation reads
\begin{equation}\label{choreographies}
q_1 (t) = q(t), q_2 (t) = q(t + T/3), q_3(t)= q(t+ 2T/3) .
\end{equation}
We call three-body solutions satisfying (\ref{choreographies}) {\it choreographies}.
The curve $q(t)$ is the curve of the eight, whose lobes
are the subject of  theorem \ref{theo:convexity}. 

 $D_6$ invariance of the figure eight implies   that  it
 is completely determined by
the three  arcs $q_1([-T/12,0]), q_2([-T/12,0]), q_3([-T/12,0])$
swept out by
the   three masses over the time  interval $[-T/12,0]$.
%
{\it In order to prove  theorem \ref{theo:convexity} it  is enough to prove that
the curvatures of these three arcs are never zero}
(with the exception of the point $q_1 (0)$ which is 
taken to be the origin
-- the self-intersection point of the eight).

\begin{figure}[htbp]
\includegraphics{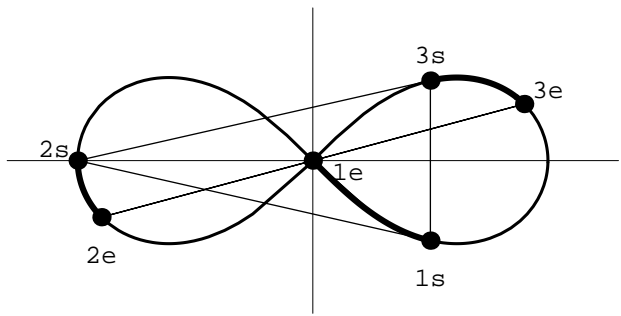}
\caption{Figure eight.
Solid circles labeled $j$s and $j$e
($ j=1,2,3$)
represent places 
at $t=-T/12$ and $t=0$
of $j$th mass.
}\label{theInterval}
\end{figure}

A configuration $(q_1, q_2, q_3)$
satisfying $q_1 + q_2 + q_3 = 0$ is called an 
{\it Euler configuration} if one of the $q_i = 0$.
Then necessarily the other two masses
$q_j, q_k$ are of the form $\zeta, -\zeta$
so that the entire configuration 
$(q_1, q_2, q_3)$
is collinear with  mass $i$
at the origin located at the midpoint  of the segment defined
by the other two masses  $j$ and $k$.  
Upon translating time if necessary, 
and relabeling mass labels, we can insist that
at the   time $0$ the configuration is an  Euler configuration
 with $1$ at the
origin and $3$ in the first quadrant as indicated in 
figure \ref{theInterval}. 
%
And  at the initial time  $t=-T/12$ the three masses form an isosceles
triangle with mass $2$ at the vertex and lying
on the negative $x$-axis. 

The eight minimizes the usual action of mechanics
(integral of the kinetic minus potential energy)
among all $T$-periodic loops enjoying $D_6$ symmetry.
Equivalently (see \cite{CM})
the path  $(q_1 (t), q_2 (t),
q_3(t))$ of the eight over
the fundamental time interval
$[-T/12,0]$   minimizes the
action among all paths starting at time
$-T/12$ in an  isosceles configuration with $2$
being the vertex and ending at 
time $0$ in an Euler configuration with $1$ being the origin.
%
An important consequence
of   minimization proved in \cite{CM}, p.896-897
is  that there are no   times in the fundamental domain
besides the endpoints 
at which the configuration is either collinear or  isosceles.
%
It follows that for all $t \in (-T/12,0)$ we have 
\begin{equation}\label{distanceOrdering}
r_{13} < r_{12} < r_{23},
\end{equation}
\begin{equation}\label{orientation}
q_1 \wedge q_2
=q_2 \wedge q_3
=q_3 \wedge q_1
<0
\end{equation}
where 
 $r_{ij} = |q_i -q_j|$ denotes the distance 
between mass $i$ and mass $j$.  
Here, we write $(x,y) \wedge
(u,v) = xv - yu$
 for planar vectors $(x,y), (u,v)$.
We call equation 
(\ref{distanceOrdering})
the  {\it distance ordering inequality}.

\subsection{Initial and Final Velocities.}\label{initialFinalVelocities}
At the Euler time, $t=0$,   the  
velocities of $2$ and $3$ are   antiparallel to the velocity of $1$, and half
its size.  See 
figure \ref{theInterval}.  
This fact follows from the action minimization
of the eight.  
%
At the isosceles time $t = -T/12$,  
$2$'s velocity  is vertical, pointing down,
and the velocities of $1$ and $3$ are such that their tangent lines
pass through $2$.  
This fact follows from the three tangents theorem
\cite{FFO}, 
and the angular momentum properties, described below.  

\subsection{Angular momentum and
star-shapedness.}\label{angularMomentumStarShapeness} 

\begin{figure}[htbp]
\includegraphics{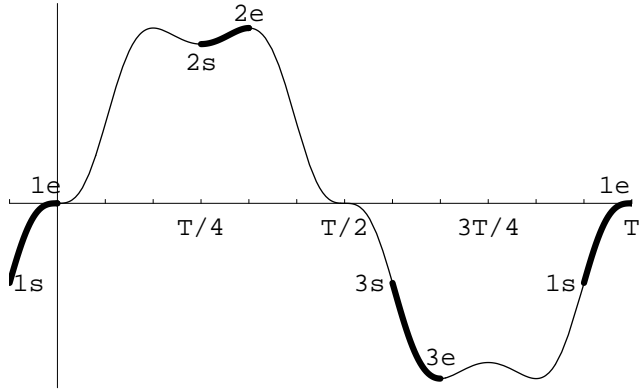}
\caption{$\ell(t)$ vs. $t$.}\label{ell}
\end{figure}

%
Write
$$\ell_j = q_j \wedge \dot q_j$$
for the angular momentum of the $j$th particle.
Action minimization of the eight implies that its total  angular momentum is zero: 
$$\ell_1 + \ell_2 + \ell_3 = 0$$
of the eight. 
Newton's equations (see \cite{CM}, p.896) imply
\begin{equation}\label{torque}
\dot \ell_3 =
\left({1 \over {r_{13} ^3 }} - {1 \over {r_{23} ^3 }}\right)(q_1 \wedge q_2)
\end{equation}
valid for all time.
Upon taking account the distance inequality
(\ref{distanceOrdering})
and (\ref{orientation})
we find that $\dot \ell_3 < 0$ on
the arc 3.
Similarly, we get
\begin{equation}\label{dotEll}
\dot \ell_1 >0, \dot \ell_2  > 0 , \dot \ell_3 < 0.
\end{equation}
By the symmetry
$\ell_{1s}=\ell_{3s}=-2\ell_{2s}<0$.
(The inequalities
$\ell_{1s}<0$ and $\ell_{1e}=0$ are consistent
with
$\dot \ell_1 >0$.)
Also
$\ell_{2s}>0$ and $\dot \ell_2  > 0$ imply
$\ell_{2e}=-\ell_{3e}>0$.
(See figure \ref{ell}.)
Therefore 
over the interior
$(-T/12,0)$ of our fundamental domain we have
\begin{equation}\label{ellInequality}
\ell_1 < 0, \ell_2 >0 , \ell_3 <0
\end{equation}
More generally,
set 
$$\ell = q \wedge \dot q$$
as $q$ varies over the eight.
It follows that on the right lobe ($x > 0$) we have
$$\ell < 0 \hbox{ for } x > 0.$$
(See figure \ref{ell}).

A curve in the plane is called `star-shaped'
with respect to the origin if every ray starting
at the origin intersects the curve at most once.
%
For a smooth curve, this is equivalent to
the assertion that,  when written in polar coordinates  
as $(r(t), \theta(t))$, the function $\theta(t)$ is strictly
monotone and does not  vary by more than $2\pi$.  
Since $\ell = r^2 \dot \theta$ the
star-shapedness of a  curve (such as one lobe
of the eight) which lies in 
the half-plane $x > 0$  is thus equivalent to $\ell \ne 0$.

\subsection{Three tangents theorem}\label{threeTangentTheorem}

The following  theorem can be found in \cite{FFO} where it was used to
find and establish existence
of a  choreographic three-body lemniscate for a non-Newtonian
potential.

\begin{theorem}[Three tangents]\label{theo:threetangent}
Let 
$(q_1 (t), q_2 (t), q_3 (t))$
be three planar
curves whose total linear and total  angular momentum 
are zero.  
Then the three instantaneous
tangent lines to these three curves
are coincident -- they all three
intersect in the same (time-dependent)
point or are parallel.
\end{theorem}

{\bf Proof. } Fix the time $t$. 
%
Because $\dot{q}_1  + \dot{q}_2 + \dot{q}_3 = 0$, 
translating all the $q_i$ in the same fixed direction
does not change the condition of having zero
angular momentum.  
So,  without loss of generality,
we can choose  the origin to be the point  of 
intersection of the tangent lines
to  $q_1$ and $q_2$ 
at   time $t$.
Because the point $q_1(t)$ lies along the line through the origin
%
%
in the direction $\dot{q}_1$ we have that
$q_1 (t) \wedge \dot{q}_{1}(t) = 0$.  
%
Similarly
$q_2 (t) \wedge \dot{q}_{2}(t) = 0$. 
%
But the
total angular momentum is zero so we must have
that $q_3 (t) \wedge \dot{q}_{3}(t) = 0$
which asserts that the line tangent to the curve
of $q_3$ at $t$ also passes through the origin.
QED

Remark:  the proof also works for unequal masses 
$m_1, m_2, m_3$.  
Simply use the correct 
mass-weighted formulae
for linear and angular momentum.

\subsection{Splitting Lemma.}\label{splittingLemma}

We will  use the following 
`Splitting Lemma' in several places in the proof. 
%
A line in the plane divides the plane
into three pieces -- two open   half-planes
and the line itself.   
We say that a point lies {\it strictly on one side}
of the line if it lies in one of the open half-planes. 
We  say that
this line {\it splits} the points 
$A$ and $B$ of the plane if the two 
points lie in opposite open half planes. 

\begin{lemma}\label{lem:splitting}
Let $(q_1 (t), q_2 (t), q_3 (t))$ 
be a planar solution to  Newton's
three-body equation with 
attractive $1/r$ potential.
Suppose that at time $t_*$
the arc $q_i (t)$ of mass $i$ has an inflection point 
and nonzero speed.  
Then
the tangent line $\ell$ to this
arc at time $t_*$ must 
either (A) split the other two
masses   $q_j (t_*)$ and $q_k (t_*)$
or  (B) all three masses must
lie on this tangent line. 
\end{lemma}

{\bf Proof.}  
Suppose, to the contrary, that either
both   $q_j (t_*)$ and $q_k (t_*)$ lie
strictly on one side of $\ell$, or that
one lies on $\ell$ while the other lies strictly on
one side.  
According to  Newton's equations  the
acceleration $\ddot q_i (t_*)$ is a linear combination of
$q_j (t_*) - q_i (t_*)$ and $q_k (t_*) - q_i (t_*)$ and  
the coefficients of this linear combination are
positive.   
Thus, translating $\ell$ and the configuration
of masses back to the origin by subtracting $q_i(t_*)$
we see that this acceleration lies strictly on one side
of the line through $0$ spanned by the velocity
$\dot q_i (t_*)$.  
Consequently,
the acceleration and velocity of $q_i (t)$ are linearly
{\it independent} at $t_*$.  
But the condition
of being an inflection point is precisely that
the acceleration and velocity be linearly dependent.
QED 

{\it Remark.}  The same proof works if the 
Newtonian potential $-\Sigma_{i < j} m_i m_j / r_{ij}$ is
replaced by any potential  $V =\Sigma_{i < j} f(r_{ij})$
where $df/dr > 0$.

\subsection{A Convexity Proposition}\label{convexityProposition}

A parameterization $t$ of a curve
$C$ is called {\it nondengenerate}
if under this parameterization the
derivative $dC(t)/dt$ is never zero. 
 A smooth, possibly self-intersecting curve
is called {\it locally convex} if its curvature never vanishes. 

\begin{proposition}\label{prop:convexity}
Let $C$ be a smooth locally convex planar curve
parameterized by a nondegenerate parameter $t$.
Let $\ell(t)$ be the 
  the tangent line to 
$C$ at $C(t)$.   Let $m$ be  a line not intersecting $C$.
Let $P(t)$
be the   point of intersection
  of $\ell(t)$ and $m$. 
Then   $dP/dt \ne 0$ for all values of $t$   (*).    
\end{proposition}

(*)  Some care is in order regarding those instants
where the line $\ell(t)$ does NOT intersect $m$ by virtue
of being parallel to it.  
The notion of convex and locally convex
is a {\bf projective} notion. 
(Local convexity is the assertion of linear 
independence of
the first and the  second derivatives of the curve,
and  this linear independence is
invariant under projective transformations.) 
So we should view the curve $C$ and
the line $m$ as lying in the real projective plane $RP_2$.  From this
new perspective,  the point $P(t)$ exists and is uniquely defined 
{\bf for all values
of} $t$.  
The theorem asserts that $dP/dt$ never vanishes.
In the proof below
we  ignore these instants of 
parallelism  i.e of $P(t)$ being at infinity
relative to the affine plane within which we are
working.  
Use  a projective transformation 
to bring the point at infinity on $m$ to a finite point on the affine plane,
and repeat the computation below to arrive at a proof valid 
for all values of $t$, including the instants of parallelism.

{\bf Proof.}  
 By a translation and rotation we can take
$m$ to be the $y$-axis.  If $(x(t), y(t))$
parameterize $C$, then the line $\ell(t)$
is given by 
$\{ (x (t), y (t)) + \lambda (\dot x (t), \dot y (t)): \lambda \in \R \}$.
It follows that   the  point of intersection 
of $\ell(t)$ with $m$ occurs at the
point $P(t) = (0, p(t)))$ where
$$p =  -\frac{x (t) \dot y (t) - y(t) \dot x (t) }{\dot x (t)}$$
%
A routine differentiation combined with the 
definition of the curvature $\kappa$ yields that
$$\frac{dp}{dt} = -{v^3 x \over \dot x ^2   } \kappa$$
where $v = \sqrt{\dot x ^2 + \dot y ^2} $
is the curve's speed.  
Every factor on the right hand side
is nonzero.   
The speed $v$ is nonzero since the curve is
nondegenerate.  
The coordinate $x$ is nonzero since
the curve $C$ never  intersects
the line $m$.  
The   only times at which 
$\dot x = 0$ are the excluded instants,
those at which $\ell(t)$ is parallel to $m$.  
And $\kappa \ne 0$
by convexity of $C$.  It follows that  $dp/dt \ne 0$ 
everywhere except at the  excluded instants.
QED

\section{To each mass its own quadrant.}\label{quadrant}

A crucial ingredient of the proof of theorem 
\ref{theo:convexity}  is
that each mass ``stays in its own quadrant''
during the time interval $(-T/12,0)$.
Initially 3 is in quadrant 1,
1 is in quadrant 4 and 2 is on the x-axis between
quadrants 2 and 3, but moving into quadrant 3.  
%
Hence, for a short time interval $(-T/12,-T/12 + \epsilon)$
mass 3 lies in quadrant 1, 1 in quadrant 4, and
2 in quadrant 3.

\begin{lemma}\label{lem:quadrant}
Over the entire time interval $(-T/12,0)$ 
body 1 lies in the 4th quadrant,
body 2 lies in the 3rd quadrant, 
and
 body 3 lies in the 1st quadrant. 
\end{lemma}

{\bf Proof of Lemma.} 

 By way of contradiction, 
suppose one of the masses leaves its initial quadrant
before the alloted time $T/12$.  
It must
exit along the boundary of this quadrant.
 It cannot exit through
the origin, as this would imply an Euler configuration
and the only Euler configuration occurs at the endpoint
of the interval.  

  We argue individually that
each mass cannot be the first exiting mass.    
%
Suppose that
$2$ exits first.  
It cannot leave crossing the x-axis as this
would contradict star-shapedness of the lobe it lies on.
Neither can it exit through the y-axis.  
For if it exited through
the y-axis, its x-coordinate would be zero.  
Now the x-coordinates
of the other two masses cannot both be zero, otherwise
this instant would be a syzygy instant.  
Thus at least one of the other
two masses lie in their quadrants, which means their x-coordinate
is positive.  
Thus the sum of the x-coordinates of the masses 
is positive, contradicting that the center of mass is at the origin.

Mass 1 cannot leave first.  
For it cannot leave through the x-axis,
as this would again contradict star-shapedness.  
It cannot leave through
the y-axis as this would violate the distance ordering 
(\ref{distanceOrdering}).  
$$r_{13} < r_{12} < r_{23}.$$  
(See figure 1.)  
To see this violation,
write the exit point for mass 1 as 
$(0,y_1)$ with $y_1 <0$.  
Then the other masses
must be at $(-x, y_2)$ and $(x, y_3)$ with
$x > 0$ (since the configuration cannot be a syzygy)
and $y_2 <0, y_3 >0$.  
We have
$r_{13}^2 = x^2 + (y_3 -y_1)^2  $,
$r_{12}^2 = x^2 + (y_2 -y_1)^2  $.
But $y_3 > 0,  0 > y_1, y_2$
and $y_1 + y_2 + y_3=0$
so that $y_3 -y_1 = -2y_1 - y_2 = 2|y_1| + |y_2|$
while $|y_2 - y_1| < |y_2| + |y_1|$ so 
that 
$(y_3 - y_1)^2 > (y_2 -y_1)^2 $
and $r_{13} > r_{12}$, contradicting the distance ordering.

%
Mass 3 cannot leave first. 
It cannot exit across the $x$-axis,
for if it did then the center of mass of the system would
have a negative $x$-coordinate, contradicting that
the center of mass is at the origin. 
It cannot leave
across the $y$-axis,
for this would contradict star-shapedness.   

Some further thought going
back through these cases  shows that we cannot have two or
more masses exiting their respective quadrants
simultaneously before the allotted time either.   
QED

\section{Proof of theorem \ref{theo:convexity}}\label{proof}

Denote   the arc swept out by mass $j$ during the
the time interval $[-T/12,0]$
as arc $j$.  
Write $\kappa_j$ for the curvature of arc $j$.
We must show that $\kappa_j \ne 0$
at each point of  each arc, 
with the exception of the origin  for arc 1.
More precisely, with our orientation and labelling
of the eight, we must show that
$\kappa_1 \le 0$ with $\kappa_1 <0$ for $t \ne 0$,
that $\kappa_2 > 0$ and that $\kappa_3 <0$.

\subsection{Convexity of arc 1}\label{arc1}
We begin by showing that
$\ddot{y}_1>0$ along arc 1.
Since each mass stays in its
own quadrant, we have 
$(y_3-y_1)>0$.
And by the distance ordering inequality 
(\ref{distanceOrdering})
 $r_{13}<r_{12}$. It follows that
\begin{eqnarray*}
 \ddot{y}_1 & = & (y_3-y_1)/r_{13}^3 +
(y_2-y_1)/r_{12}^3
\\
 & > &(y_3-y_1)/r_{12}^3+(y_2-y_1)/r_{12}^3 \\
 & = & -3y_1/r_{12}^3 \cr
& > & 0.
\end{eqnarray*}

Next we show that 
$\dot y_1 > 0$
along the arc.  
From the fact that
$\ddot y_1 >0$, it suffices 
to show that $\dot y_1 >0$ at
the initial point of arc 1, the isosceles point.
By the three tangents theorem and the fact that
$\ell_1 <0$ it follows that 
at the isosceles point
$\dot q_1$ points from
$q_1$ to the vertex $q_2$, so that 
$\dot y_1 >0$.

We have seen above that  $\ell_1 <0$
while 
$\dot{\ell}_1>0$ along the arc.
Combining these inequalities, we see that 
$\dot{\ell}_{1}\dot{y}_{1}-\ell_{1}\ddot{y}_{1}>0$
holds along the arc.
 On the other hand,
expanding the angular momentum,
we get
$\dot{\ell}_{1}\dot{y}_{1}-\ell_{1}\ddot{y}_{1}
=(x_{1}\ddot{y}_{1}-y_{1}\ddot{x}_{1})\dot{y}_{1}-(x_{1}\dot{y}_{1}-y_{1}\dot{x}_{1})\ddot{y}_{1}
=y_{1}(\dot{x}_{1}\ddot{y}_{1}-\dot{y}_{1}\ddot{x}_{1})
=y_{1} v_{1}^{3} \kappa
$.
Thus $y_{1} v_{1}^{3} \kappa_{1}>0$.
Since $y_1 < 0, v_1 >0$ we have  that   $\kappa_{1}<0$.

\subsection{Convexity of arc 2}\label{arc2}
Assume, by way of contradiction,  that there exists an
 inflection point  $\kappa_2=0$
on   arc 2. 
Let $a$ be the last inflection point  on arc 2 -- the one  
whose time $t$ is
closest to $0$. 
%
From the initial conditions at $t=-T/12,0$ described above
we also know that
$\kappa_2>0$ at the points $2s$ and $2e$.
By continuity, $\kappa_2 >0$ near both of these points.
Then $\kappa_2 >0$ on the arc $a \rightarrow 2e$.

\begin{figure}[htbp]
\includegraphics{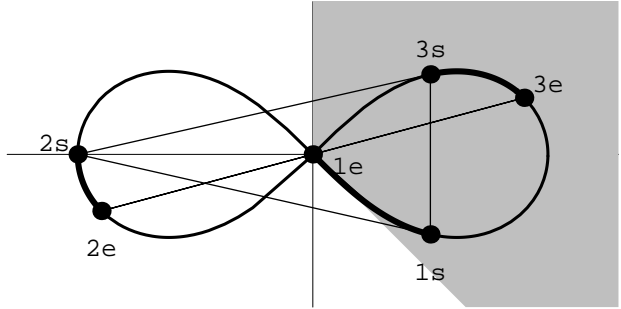}
\caption{Region for bodies 1 and 3.}\label{bodies13}
\end{figure}

We  know by the previous subsection that
 arc 1 is convex ($\kappa_1 <0$)
and  we also know that body 3 moves in 
`its own quadrant' - the 1st quadrant.
It follows  that  bodies 1 and 3 must lie within the shaded region
in the figure~\ref{bodies13}.
%

Consider the Gauss map (hodograph) of 
arc 2.  
This is the map which
assigns to a point of arc 2 the unit
tangent to arc 2,
$\dot{q}_{2}/|\dot{q}_{2}|$,
at that point.
%
\begin{figure}[htbp]
\includegraphics{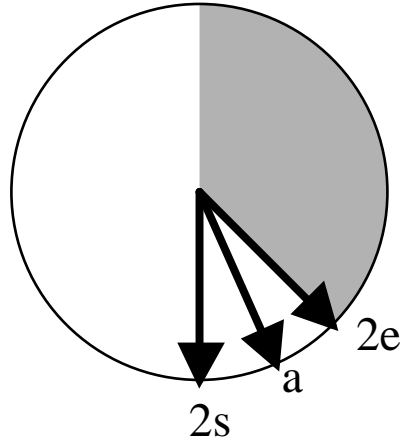}
\caption{Gauss map of 
the unit tangent vector
$\dot{q}_{2}/|\dot{q}_{2}|$.}\label{gaussMap}
\end{figure}

By Newton's equation  and the fact that
$x_1- x_2$ and $x_3 - x_2$ are positive
we have that $\ddot x_2 > 0$ on the entire arc 2.
Since $\dot x_2 = 0$ at $2s$, this implies
that $\dot x_2 >0$ on the open arc of 2, from $2s$ to $2e$,
and so in particular $\dot x_2 > 0$ at  $a$. 
Since $\kappa_2 > 0$
on the arc $a \to 2e$, the vector $\dot{q}_{2}/|\dot{q}_{2}|$
must approach  $2e$ from the point $a$
monotonically counterclockwise.  
Therefore the point $a$
lies on the arc between the points $2s$ and $2e$ on the
right half of the circle as shown in the Gauss map
(figure \ref{gaussMap}). 

But then the tangent line to arc 2 at $a$
cannot split the points $1$ and $3$,
which, according to the splitting lemma (sec. \ref{splittingLemma}), 
contradicts the assumption that
$a$ is an inflection point. 

Thus we have proved that there is no inflection point on the arc 2.
In other word, $\kappa_2 >0$ on the arc 2.

\subsection{Convexity of arc 3}\label{arc3}
Assume, by way
of contradiction,  that there are one
or more inflection points  $\kappa_3=0$
on the arc 3. 
Let $b$ be the first such point, 
the one for which the time $t$ is closest to $-T/12$. 
%
Then, by the splitting lemma (sec. \ref{splittingLemma}),  the tangent line
to arc 3 at b must split bodies 1 and 2.
In order to do that, the line  must have 
passed earlier through either body 1 or body 2.
We argue that both passings are  
impossible.

The tangent line to arc 3 cannot pass through body 1.
For, by the three tangent theorem, at the instant this
happened, 
the tangent line from the body 2
would also  pass through the body 1.
We have already proved that $\kappa_2>0$ on the arc 2.
Thus the tangent line from the body 2
never pass through the body 1 in this interval.
(See figures \ref{bodies13} and \ref{gaussMap}.)
This is a contradiction.

The tangent line to
arc 3 cannot pass through body 2.  
For if it did,  by the
three tangents theorem, the tangent line to 1's curve
would also pass through body 2 at the same instant. 
To
see that this latter passing is impossible,  join the endpoints
2s and 2e of arc 2 by a straight line $m$.  
Arc 2
lies completely on one side of this line, by
convexity.  

\begin{figure}[htbp]
\includegraphics{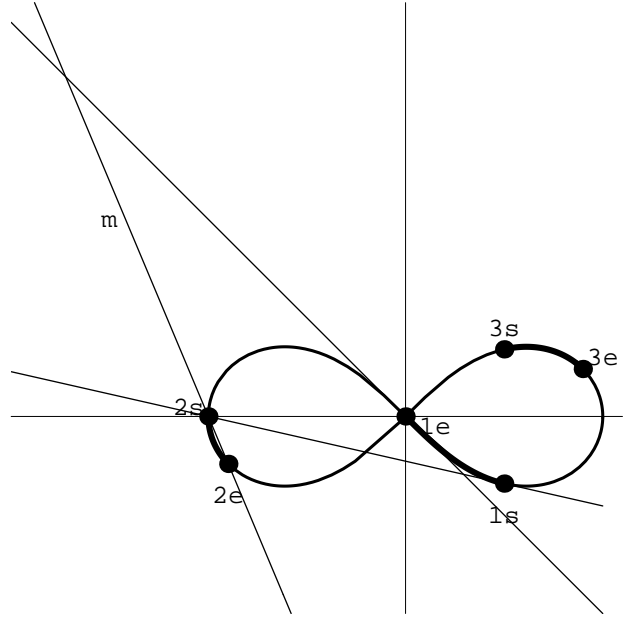}
\caption{Line $m$ and tangent lines to arc 1
at $t=-T/12$ and $t=0$.}\label{MandC}
\end{figure}
In section \ref{convexityProposition}
we proved the proposition
that if
 $c(t)$ is a smooth convex curve parameterized
at a nonzero speed, and if $m$
is a fixed line, then the intersection point $P(t)$ of $c$'s tangent line
with $m$ moves monotonically: $dP/dt \ne 0$. 
We apply the proposition to our situation. 
  At the final points e,
the tangents to 1 and 2 are parallel, so that the
intersection of $m$ with $1$'s tangent lies in 
the massless quadrant $x< 0, y > 0$.
At the initial point $s$ the intersection point
of $m$ and arc $1$'s tangent is $2s$.  
Consequently,
in between $s$ and $e$ the intersection always
lies in that part of $m$ lying in the massless quadrant.
But in order for $1$'s tangent to 
pass through 2, $1$'s tangent would have to cross
line $m$ between $2s$ and $2e$, which is in the quadrant
of arc 2, and hence it is impossible that
this tangent passes through $2$.

Therefore, we have proved that 
there is no inflection point on the arc 3.
In other word,
$\kappa_3<0$ on the arc 3.

\subsection{Conclusion}\label{conclusion}
Combining 4.1, 4.2, and 4.3 proves theorem 1.

\section{Convexity for other potentials}\label{otherPotential}
Theorem \ref{theo:convexity} holds for the ``eights'' of other potentials.
Indeed, our proof only depended on the properties and propositions
of the eight listed in section 2 
and a ``monotonicity'' property of the Newtonian potential
discussed below. 

To be precise, we  need to   define what we mean by an ``eight''
Let 
$$V = V(r_{12}, r_{23} , r_{31}) $$ 
be a three-body potential  depending
only on the interparticle distances $r_{ij}$ and invariant
under interchange of the masses.  
Then the symmetry group $D_6$ of the
eight acts on solutions to the corresponding Newton equation, taking solutions
to solutions, and so we can
speak of {\it $D_6$-invariant solutions}.

A planar solution to the Newton's equation
for $V$  will be called an {\it eight solution}
if

i) it is invariant under the $D_6$ symmetries 

ii)  on 
the interior of each fundamental
domain ($(mT/12, (m+1) T/12)$,
$m = 0,\pm 1, \pm 2, ...$) the configuration is
never  collinear and never 
isosceles 

iii) the solution  has no collisions 

Such a solution will necessarily be a planar choreography
(see 
(\ref{choreographies}) 
above), and so the three masses travel a single planar curve.
 Condition (i) implies that the center of mass is $0$
and that the angular momentum is zero.  
If, in addition,
our potential $V$ 
has the form
$$ V =\sum_{i<j} f(r_{ij})$$
with 

iv)
$df/dr > 0$ (attractive two-body potential) 
and  with
 
v)
$g(r) = r^{-1} df/dr$ 
a strictly monotone decreasing function of $r$, \\
then
all properties and inequalities 
used in this paper
hold. 
 
Indeed, return to   the starting point,
the distance ordering inequality (\ref{distanceOrdering}).
%
At $t=-T/12$ and $t=0$ we have
$r_{2s3s}=r_{1s2s}$ and 
$r_{1e2e}=r_{3e1e}<r_{2e3e}(=2r_{1e2e})$.
By the property (ii), 
the possible distance orderings on
the time interval $(-T/12,0)$ are
$r_{31}<r_{12}<r_{23}$ or
$r_{12}<r_{31}<r_{23}$.
Consider the equation for $\dot{\ell}_1$,
$$\dot{\ell}_1=(g(r_{21})-g(r_{31}))(q_2 \wedge q_3)$$
for a  monotone decreasing function $g(r)$.
We have $\dot{\ell}_1>0$ for the first ordering
and $\dot{\ell}_1<0$ for the second ordering.
But, since
$\ell_{1s}<0$ and $\ell_{1e}=0$,
$\dot{\ell}_{1}$ must be positive.
So we must have the first ordering,
namely,
the equation~(\ref{distanceOrdering}).
%
Then,
all equalities and inequalities in this paper
hold upon   replacing 
$1/r^{3}$ with $g(r)$.
Thus we have the following theorem.

\begin{theorem}\label{theo:convexityVf}
Let $V$ be a three-body potential
of the form  $V =\sum_{i<j} f(r_{ij})$ where $f$ satisfies
 (iv) and (v) immediately above, and
 admitting an eight solution as defined
by (i)-(iii) above. 
 Then each lobe of this eight for $V$  is convex.
\end{theorem}

The theorem begs the question,
do eight solutions exist for 
any potentials besides Newton?
Recall from \cite{CM}, p. 896-897
that if  a solution which satisfies (i) and (ii)  is known
to minimize the action associated to $V$ among all paths
satisfying (i), and if that solution is not identically collinear,
then automatically the solution satisfies (ii).  
The power law potentials
$$V_a = (a)^{-1}(r_{12}^a + r_{23}^a + r_{31}^a), $$
 for $a \le -2$
admit such collision-free action minimizing solutions,
and consequently they admit eight solutions.
%
Moreover, the proof of \cite{CM}, specific to   $a = -1$,
is based on
strict inequalities, and hence is
valid for a range of exponents
$-1 - \epsilon_1 < a < -1 + \epsilon_2$
for $\epsilon_1, \epsilon_2$ positive numbers. 

As a corollary, we obtain
\begin{corollary} 
For the power law potentials
$V_a$ with $a < -2$ or
with $a$ in some open interval about $-1$,
there exist eight solutions and each lobe
of these eight solutions is convex. 
\end{corollary}

\end{document}